\documentclass[12pt,leqno]{article}
\usepackage{amsmath, amsthm, amsfonts, amssymb, color}
\usepackage{mathrsfs}
\usepackage[notcite,notref]{showkeys}

\newtheorem{thm}{Theorem}[section]

\newtheorem{lem}[thm]{Lemma}

\newtheorem{Remark}[thm]{Remark}
\theoremstyle{definition}

\newcommand{\scr}[1]{\mathscr #1}
\definecolor{wco}{rgb}{0.5,0.2,0.3}

\numberwithin{equation}{section}

\newcommand{\ua}{\uparrow}

\title{{\bf 
Approximation to Wiener measure on a general noncompact Riemannian manifold}\footnote{Supported in
 part by  NNSFC (11371099).} }
\author{
{\bf   Bo Wu}\\
\footnotesize { School  of Mathematical Sciences, Fudan
University, Shanghai 200433, China}\\
 \footnotesize{Institute for Applied Mathematics, University of Bonn,
Endenicher Allee 60, 53115 Bonn, Germany}\\
\footnotesize{wubo@fudan.edu.cn, wu@iam.uni-bonn.de}}

\date{}
\begin{document}
\maketitle
\def\paral{/\kern-0.55ex/}
\def\parals_#1{/\kern-0.55ex/_{\!#1}}
\def\R{\mathbb R} \def\EE{\mathbb E} \def\N{\mathbb N} \def\Z{\mathbb Z} \def\ff{\frac} \def\ss{\sqrt}
\def\H{\mathbb H}
\def\dd{\delta} \def\DD{\Delta} \def\vv{\varepsilon} \def\rr{\rho}
\def\<{\langle} \def\>{\rangle} \def\GG{\Gamma} \def\gg{\gamma}
\def\ll{\lambda} \def\LL{\Lambda} \def\nn{\nabla} \def\pp{\mathbb P} \def\ss{\mathbb S}
\def\d{\text{\rm{d}}} \def\Id{\text{\rm{Id}}}\def\loc{\text{\rm{loc}}} \def\bb{\beta} \def\aa{\alpha} \def\D{\scr D}
\def\E{\scr E} \def\si{\sigma} \def\ess{\text{\rm{ess}}}
\def\beg{\begin} \def\beq{\beg}  \def\F{\scr F}
\def\Ric{\text{\rm{Ric}}}
\def\Var{\text{\rm{Var}}}
\def\Vol{\text{\rm{Vol}}}
\def\Scal{\text{\rm{\bf Scal}}}
\def\Ent{\text{\rm{Ent}}}
\def\Hess{\text{\rm{Hess}}}\def\B{\scr B}
\def\e{\text{\rm{e}}} \def\ua{\underline a} \def\OO{\Omega} \def\b{\mathbf b}
\def\oo{\omega}     \def\tt{\tilde} \def\Ric{\text{\rm{Ric}}}
\def\cut{\text{\rm{cut}}} \def\P{\mathbb P} \def\ifn{I_n(f^{\bigotimes n})}
\def\fff{f(x_1)\dots f(x_n)} \def\ifm{I_m(g^{\bigotimes m})} \def\ee{\varepsilon}
\def\C{\scr C}
\def\M{\scr M}\def\ll{\lambda}
\def\X{\scr X}
\def\P{\scr P}
\def\T{\scr T}
\def\A{\mathbf A}
\def\LL{\scr L}\def\LLL{\Lambda}
\def\gap{\mathbf{gap}}
\def\div{\text{\rm div}}
\def\Lip{\text{\rm Lip}}
\def\dist{\text{\rm dist}}
\def\cut{\text{\rm cut}}
\def\supp{\text{\rm supp}}
\def\Cov{\text{\rm Cov}}
\def\Dom{\text{\rm Dom}}
\def\Cap{\text{\rm Cap}}\def\II{{\mathbb I}}\def\beq{\beg{equation}}
\def\sect{\text{\rm sect}}\def\H{\mathbb H}

\begin{abstract}

In prior work \cite{AD} of Lars Andersson and Bruce K. Driver, the path space with finite interval over a compact Riemannian manifold is approximated by finite dimensional manifolds $H_{x,\P} (M)$ consisting of piecewise geodesic paths adapted
to partitions $\P$ of $[0,T]$, and the associated Wiener measure is also approximated by a sequence of probability measures on finite dimensional manifolds. In this article, we will extend their results to the general path space(possibly with infinite interval) over a non-compact Riemannian manifold by using the cutoff method of compact Riemannian manifolds. Extension to the free path space. As applications, we obtain integration by parts formulas in the path space $W^T_x(M)$ and the free path space $W^T(M)$ respectively.
\end{abstract}

\noindent Keywords: Riemannian manifold; Wiener measure; Brownian motion;
Path space.\vskip 2cm

\section{Introduction}\label{sect1}
Let $(M,g)$ be an $n$-dimension complete Riemannian manifold and $d$ be the Riemannian distance on $M$. For each fixed $T>0$ and $x\in M$,  the based path space and the path space with a fixed point over $M$ are given by
$$\aligned &W^T(M):=C([0,T]: M),
\\& W^T_x(M):=\big\{\gamma\in W^T(M): \gamma_0=x\big\}.\endaligned $$
Then $W^T(M)$ and $W^T_x(M)$ are respectively Polish
spaces under the uniform distance
$$\rho(\gamma,\sigma):=\displaystyle\sup_{t\in[0,T]}d(\gamma_t,\sigma_t),\quad\gamma,\sigma\in W^T(M).$$
In particular, let $\rho_x(\gamma):=\rho(\gamma,x), \gamma\in W^T_x(M)$ be the distance function on $W^T_x(M)$ starting from $x$.

Denote by $O_x(M)$ be the orthonormal frame bundle at $x\in M$, then $O(M):=\sup_{x\in M}O_x(M)$ is the orthonormal frame bundle over $M$. Let $U_t^x$ be the horizontal Brownian motion on $O(M)$ generated by the laplace operator $\Delta$; that is, $U_t^x$ solves
the following stochastic differential equation on $O(M)$,
\begin{equation}\label{eq1.1}
\d U_t^x=\sqrt{2}\displaystyle\sum^d_{i=1}H_i(U_t^x)\circ\d W_t^i, \quad U_0\in O_x(M),
\end{equation}
where $W_t=(W_t^1,\cdots,W_t^d)$ is the $d$-dimensional Brownian motion, $\{H_i\}_{i=1}^n: TM\rightarrow TO(M)$ is a standard othornormal basis
of horizontal vector fields on $O(M)$. Let $\pi: O(M) \rightarrow M$ be the canonical projection.
Then $X^x_t:=\pi(U^x_t),\ t\geq0$ is the Brownian motion with initial point $x$. Throughout the article, we assume that $M$ is complete and stochastic complete, which implies that $X^x_t$ is non-explosive. 

Let us define $\F C_0^T$ the class of cylindric smooth functions on the path space $W^T(M)$ defined as
\begin{equation}\label{eq1.2}\aligned \F C_0^T:=\Big\{F(\gamma)&=f(\gamma(t_1),\cdots,\gamma(t_N)):\ N\geq1,~\gamma\in W^T(M),\\&
~~~~~~~~~~~~~~~~~~~
0<t_1<t_2\cdots<t_N\leq T,f\in C_0^{Lip}(M^N)\Big\},\endaligned\end{equation}
where $C_0^{Lip}(M^N)$ is the set of all Lipschitz continuous functions with compact supports on $M^N$. Denote by $\pp_x^T$ the Wiener measure on $W^T_x(M)$, that is,
 \begin{equation}\label{eq1.3}\aligned \int_{W^T_x(M)}F(\gamma)\d \pp_x^T(\gamma)=\int_{M^N}f(y_1,\cdots,y_N)\prod_i^Np_{\Delta_i}(x_{i-1},x_i)\d x_1\cdots\d x_N,\endaligned\end{equation}
for any $F\in\F C_b^T$ of the form $F(\gamma)=f(\gamma(t_1),\cdots,\gamma(t_N))$, where $p_t(x,y)$ is the heat kernel on $M$ and $\Delta_i=t_i-t_{i-1}$.

Before stating our main result, we need to introduce some notation partly following from \cite{AD}. Let $\T^T$ be the set of all partitions of $[0,T]$ and $E(\gamma):=\int^T_0\langle\gamma'(t),\gamma'(t)\rangle\d t$ be the energy of each absolute  continuous path $\gamma\in W^T_x(M)$, otherwise let $E(\gamma)=+\infty$.
The path space with finite energy is define by
\begin{equation}\label{eq1.4}\aligned
H_x^T(M):=\left\{\gamma\in W^T_x(M): E(\gamma)<\infty\right\}.
\endaligned\end{equation}
For every $\P\in \T^T$, define
\begin{equation}\label{eq1.5}\aligned
H^T_{x,\P}(M):=\left\{\gamma\in H^T_x(M)\cap C^2([0,T]/ \P): \nabla \gamma'(t)/ \d t=0, s\notin \P\right\}.
\endaligned\end{equation}
The tangent space of $H^T_x(M)$ at $\gamma$ may be naturally identified with the space of all absolutely continuous vector fields $X:[0,T]\rightarrow TM$ along $\gamma$ such that $X(0)=0$ and 
\begin{equation}\label{eq1.6}\aligned
G_x^{T,1}(X,X):=\int^T_0\left\langle\frac{\nabla X(t)}{\d t},\frac{\nabla X(t)}{\d t}\right\rangle\d t,
\endaligned\end{equation}
where
\begin{equation}\label{eq1.7}\aligned
\frac{\nabla X(t)}{\d t}:=\parals_t(\gamma)\frac{\d}{\d t}(\parals_t^{-1}(\gamma)X(t)).
\endaligned\end{equation}
Here $\parals_t(\gamma):T_xM\rightarrow T_{\gamma(t)}M$ is parallel translation along $\gamma$ relative to Levi-Civita covariant derivative $\nabla$. Similarly, we may define a weak Riemannian metric $G_x^{T,0}$ on $H^T_x(M)$ by
\begin{equation}\label{eq1.8}\aligned
G_x^{T,0}(X,X):=\int^T_0\left\langle X(t),X(t)\right\rangle \d t,
\endaligned\end{equation}
for any $X\in TH^T_x(M)$. By the argument in \cite{AD}, these metrics may be interpreted $D\gamma$ as the `Riemannian volume measures' $\Vol_{G_x^{T,1}}$ and $\Vol_{G_x^{T,0}}$ with respect to $G_x^{T,1}$ and $G_x^{T,0}$ respectively.
By induction, we may get the metrics on $H^T_{x,\P}(M)$ for some partition $\T^T\ni\P=\{0=t_0<t_1<\cdots<t_N=T\}$,
\begin{equation}\label{eq1.9}\aligned
G^{T,1}_{x,\P}(X,Y):=\sum^N_{i=1}\left\langle\frac{\nabla X(t_{i-1}+)}{\d t},\frac{\nabla Y(t_{i-1}+)}{\d t}\right\rangle\Delta_it,
\endaligned\end{equation}
and 
\begin{equation}\label{eq1.10}\aligned
G^{T,0}_{x,\P}(X,Y):=\sum^N_{i=1}\left\langle X(t_{i-1}+),Y(t_{i-1}+)\right\rangle\Delta_it,
\endaligned\end{equation}
for all $X,Y\in TH^T_{x,\P}(M)$.

Assume $\Vol_{G^{T,1}_{x,\P}}$ and $\Vol_{G^{T,0}_{x,\P}}$ are respectively the volume forms on $H^T_{x,\P}(M)$ with respect to $G_{x,\P}^{T,1}$ and $G_{x,\P}^{T,0}$ and the associated probability measures are defined by
\begin{equation}\label{eq1.11}\aligned
&\pp_{x,\P}^{T,1}:=\frac{1}{Z^1_{\P,T}}\e^{-E/2}\Vol_{G^{T,1}_{x,\P}}\\
&\pp_{x,\P}^{T,1}:=\frac{1}{Z^0_{\P,T}}\e^{-E/2}\Vol_{G^{T,0}_{x,\P}},
\endaligned\end{equation}
where $Z^0_{\P,T}$ and $Z^1_{\P,T}$ are normalization constants given by
\begin{equation}\label{eq1.12}\aligned
&Z^0_{\P,T}:=\prod^N_{i=1}(\sqrt{2\pi}\Delta_it)^n\\
&Z^1_{\P,T}:=(\sqrt{2\pi})^{nN}.
\endaligned\end{equation}
For any partition $\P=\{t_1,\cdots,t_N\}$ of $\T^T$, denote by 
$|\P|:=\max\{|\Delta_i t|:i=1,\cdots,N\}$ the mesh size of the partition $\P$.

In prior of these preparations, we may now state our main results of this paper.

\beg{thm}\label{T1.1} Assume that $M$ is complete and stochastic complete. Then 

$(1)$ For each bounded continuous function $F$, we have
\begin{equation}\label{eq1.13}\aligned
\lim_{|\P|\rightarrow0}\int_{H^T_{x,\P}(M)}F(\gamma)\d \pp^{T,1}_{x,\P}(\gamma)
=\int_{W^T_x(M)}F(\gamma)\d \pp_x^T(\gamma)
\endaligned\end{equation}

$(2)$ 
Suppose that 
$$\int_{W^T_x(M)}\e^{-\frac{1}{6}\int^T_0\Scal(\gamma(t))\d t}\d \pp_x^T(\gamma)<\infty.$$
Then for  each bounded continuous function $F$, we have
\begin{equation}\label{eq1.14}\aligned
\lim_{|\P|\rightarrow0}\int_{H^T_{x,\P}(M)}F(\gamma)\d \pp^{T,0}_{x,\P}(\gamma)
=\int_{W^T_x(M)}F(\gamma)\e^{-\frac{1}{6}\int^T_0\Scal(\gamma(t))\d t}\d \pp_x^T(\gamma),
\endaligned\end{equation}
where $\Scal$ is the scalar curvature of $M$.
\end{thm}

\begin{Remark}\label{r1.2} $(a)$ When $M$ is compact and $T=1$, it is proven in \cite{AD}.

$(b)$ In fact, it is not difficult to see that Theorem 1.8 
in \cite{AD} still holds when the time $1$ is replaced by a constant $T>0$.
\end{Remark}

Since Feynman\cite{F} established the path integral formula which play a crucial role of the theory in quantum mechanics, much of the current interest concerning the direction has grown deeply. The role of path integrals in quantum mechanics is surveyed by Gross in [23] and detailed more by Feynman\cite{F}as well as Glimm-Jaffe\cite{GJ}.  When $M=\R^d$, it is known that Wiener measure on $W_x(\R^d)$  may be approximated by Gaussian measures on piecewise linear path spaces. When $M$ is compact and $T=1$, Anderson-Driver in \cite{AD} gave two different approximations of Wiener measure. These results  gives a rigorous interpretation of a Feynman path integral on a non-compact Riemannian manifold $M$, which is heuristically of the form
$$\frac{1}{Z^0}\int_{W^T_x(M)}F(\gamma)\e^{-\frac{1}{6}\int^T_0\Scal(\gamma(t))\d t}\d \D$$
and
$$\frac{1}{Z^1}\int_{W^T_x(M)}F(\gamma)\e^{-E(\gamma)/2-\frac{1}{6}\int^T_0\Scal(\gamma(t))\d t}\d \D,$$
where $Z^0,Z^1$ is ``normalization constants", the reader may refer to \cite{AD} for the detailed argument. After that, 
their result are extended to compact manifold or other manifolds with certain curvature conditions by \cite{TL,AL,LZH}.  

The rest of this paper is organized as follows: In Section 2, we will present the
proof of Theorem \eqref{T1.1}, and by which we will get the formula of integration by parts on path space. Finally, extension for the infinite interval will provided in Section 3.

\section{Approximation to Wiener measure for finite interval}
\subsection{Proof of Theorem \ref{T1.1}}
In this subsection, we will prove Theorem \ref{T1.1}. To do that, we first introduce a Lemma proven in \cite{CLW1} about cutoff of compact Riemannian manifold.  Since $M$ is complete,
there exists a $C^\infty$ non-negative smooth function $\hat d:M\rightarrow\R$ with the property that $0<|\nabla \hat d|\le 1$  and
$$\left|\hat d(x)-\frac{1}{2}d(x)\right|<1,\quad \forall\; x\in M.$$
For every non-negative $m$, define $D_m:=\hat d^{-1}((-\infty,m)):=\{z \in M; \hat d(z)<m\}$, then it is easy to verify
$B_x(2m-2)\subset D_m\subset  B_x(2m+2)$, where $B_x(r):=\{z \in M; d(z)<r\}$ is the geodesic ball centred at $x$ of radius $r$.

\begin{lem}\label{l2.1}[Lemma 5.1 in \cite{CLW1}]
For every $m \in \mathbb{Z}_+$, there exists a smooth compact Riemannian manifold $(\tilde M_m,\tilde g_m)$, such that 
$(D_m,g)$ is 
isometrically embedded into $(\tilde M_m,\tilde g_m)$ as an open set. 
\end{lem}

\ \newline\emph{{\bf Proof of Theorem \ref{T1.1}.}}  
Here we only prove $(1)$($(2)$  may be dealt similarly). In order to do that, we will take a sequence of compact Riemannian manifolds  $\{M_m\}_{m\geq1}$ such that $M_m$ converges to $M$ under suitable sense. In addition, combining this with Anderson-Driver's results for compact Riemannian manifold, we will obtain the aim. 
By lemma 2.1, we know that there exists a sequence of  compact Riemannian manifolds $\{(M_m, g_m)\}_{m\geq1}$ such that for each $m\geq1$, $D_m$ is isometrically embedded into some open set of $(M_m, g_m)$. 

In order to compare the Browinian motions in different Riemannian manifolds $M_m$, we need to model them into the same probability space. Suppose that
$(\Omega, \F, \P)$ is a complete probability space, and $W_t$ is a $\R^n$-valued Brownian motion on this space.
We consider the SDE on the frame bundle space $O(M)$,
\begin{equation}\label{eq2.1}
\begin{cases}
&\d U_t=\displaystyle\sum^n_{i=1}H_i(U_t)\circ\d W_t^i,\ \ t \in [0,T],\\
& U_0=u_x,
\end{cases}
\end{equation}
where $u_x\in O_x(M)$ with $\pi(u_x)=x$. Let $X_t:=\pi(U_t)$, then  $X_t$ is the Brownian motion on $M$ starting from $x$, and $U_{\cdot}$ is the horizontal lift along
$X_{\cdot}$.
Under the sense of isometry, we may look $B_x(m)$ as subset of $M$ and $M_m$. 
Similarly, since $\langle\cdot, \cdot\rangle_m= \langle\cdot, \cdot\rangle$ on $B_x(m)$,
we can choose an orthonormal basis $\{H_{i,m}\}_{i=1}^n$  of horizontal vector fields
on $O(M_m)$ such that $H_{i,m}(u)=H_{i,R}(u)=H_i(u)$ for every $R\ge m$ when $u \in O(M_m)$ satisfies
$\pi(u) \in B_x(m)$. Let $W_t$, $u_0$ be the same as that in (\ref{eq2.1}),
we consider the following SDE,
\begin{equation*}
\begin{cases}
&\d U_{t,m}=\displaystyle\sum^n_{i=1}H_{i,m}(U_{t,m})\circ\d W_t^i,\ \ t \in [0,T],\\
& U_{0,R}=u_o,
\end{cases}
\end{equation*}
so $X_{\cdot,m}:=\pi(U_{\cdot,m})$ is the Brownian motion on $M_m$, $U_{\cdot,m}$ is the horizontal lift
along $X_{\cdot,m}$ on $M_m$. Moreover, let $\tau_m:=
\inf\{t\ge 0: X_t \notin B_x(m)\}$,
we have $U_{t,m}=U_{t,m}=U_t$ $\P$-$a.s.$ for every
$R \ge m, t \le \tau_m$. Let $\pp_x^{T,m}$ be the distribution of the Brownian motion $X_{\cdot,m}$ on $M_m$. Then for any bounded continuous function $F$ with $F(\gamma)=0$ if $\gamma\cap B_x^c(m)\neq \emptyset$, we have
\begin{equation}\label{eq2.2}\aligned
\pp_x^T(F)=\pp_x^{T,m}(F).
\endaligned\end{equation}

For any $\P\in \T^T$, denote respectively by $W^T_x(M_m), H^T_{x,\P}(M_m),\pp^{T,0,m}_{x,\P}$ and $\pp^{T,1,m}_{x,\P}$ with $M$ replaced by $M_m$.
To prove $\eqref{eq2.13}$, we only show that $\eqref{eq2.13}$ holds for any bounded non-negative continuous function $F$. In fact, since $F=F^+-F^-$, we obtain easily our conclusion. For every $m\ge 2$, choose a
$l_m\in C^\infty_0(\R)$ such that
\begin{equation}\label{eq2.3}
l_m(r)=\begin{cases}
& 1,\quad\quad\quad\quad\quad\text{if}\ |r|\le m-1,\\
& \in [0,1], ~\quad\quad~\text{if}\ m-1<|r|<m,\\
&0,\quad\quad\quad\quad\quad\text{if}\ |r|\ge m.
\end{cases}
\end{equation}
Let
$d_m$ be the Riemannian distance on $M_m$ and
$$\rho_m(\gamma):=\sup_{t \in [0,T]}d_m(\gamma(t),o), \ \ \phi_m(\gamma):=l_m(\rho_m(\gamma)),\quad \gamma \in W^T_x(M_m).$$
Note that $l_m(r)=0$ if $r \ge m$ and $\rho_m(\gamma)=\rho(\gamma)$
for each $\gamma\in W_x(M_m)\cap W^T_x(M)$ with $\gamma\subset B_x(m)$,  so we can extend $\phi_m$ to be defined in
$W^T_x(M)$ by $\phi_m(\gamma):=l_m(\rho(\gamma))$ for $\pp_x^T$-$a.s.$ $\gamma \in W^T_x(M)$.
Let $F_m:=\phi_m F$, then we obtain
\begin{equation}\label{eq2.4}\aligned
\lim_{|\P|\rightarrow0}\int_{H^T_{x,\P}(M)}F_R(\gamma)\d \pp^{T,1}_{x,\P}(\gamma)&=\lim_{|\P|\rightarrow0}\int_{H^T_{x,\P}(M_m)}F_m(\gamma)\d \pp^{T,1,m}_{x,\P}(\gamma)
\\&=\int_{W^T_x(M_m)}F_k(\gamma)\d \pp_x^{T,m}(\gamma)=\int_{W^T_x(M)}F_R(\gamma)\d \pp_x^T(\gamma),
\endaligned\end{equation}
where the first inequality is due to the same underlying metric on $B_x(m)$ and the cutoff function $\phi_m$; the second one follows from Anderson-Driver's results; the final inequality apply \eqref{eq2.2}. By using the triangle inequality, we have
\begin{equation}\label{eq2.5}\aligned
&\left|\int_{H^T_{x,\P}(M)}F(\gamma)\d \pp^{T,1}_{x,\P}(\gamma)-\int_{W^T_x(M)}F(\gamma)\d \pp_x^T(\gamma)\right|
\\\leq 
&\left|\int_{H^T_{x,\P}(M)}F(\gamma)\d\pp^{T,1}_{x,\P}(\gamma)
-\int_{H^T_{x,\P}(M)}F_m(\gamma)\d \pp^{T,1}_{x,\P}(\gamma)\right|
\\&+\left|\int_{W^T_x(M)}F_m(\gamma)\d \pp_x^T(\gamma)
-\int_{W^T_x(M)}F(\gamma)\d \pp_x^T(\gamma)\right|\\
&+\left|\int_{H^T_{x,\P}(M)}F_m(\gamma)\d \pp^{T,1}_{x,\P}(\gamma)
-\int_{W^T_x(M)}F_m(\gamma)\d \pp_x^T(\gamma)\right|
\\\leq& 
\int_{H^T_{x,\P}(M)}(1-\phi_m)F(\gamma)\d\pp^{T,1}_{x,\P}(\gamma)
+\int_{W^T_x(M)}\left|F-F_m\right|(\gamma)\d \pp_x^T(\gamma)\\
&+\left|\int_{H^T_{x,\P}(M)}F_m(\gamma)\d \pp^{T,1}_{x,\P}(\gamma)
-\int_{W^T_x(M)}F_m(\gamma)\d \pp_x^T(\gamma)\right|\\=:&I^\P_1(m)+I_2(m)
+I^\P_3(m).
\endaligned\end{equation}
Since $M$ is stochastic complete, 
$$\lim_{R\rightarrow \infty}\pp_x^T(\rho_x\geq m)=0.$$ 
Combining this with
$$I_2(R)\leq \Vert F\Vert_\infty\pp_x^T(\rho_x\geq m),$$
for any $\varepsilon>0$, there exists a large enough constant $m_0$ such that 
$I_2(m_0)\leq \varepsilon/3$ and
\begin{equation}\label{eq2.6}\aligned
\pp_x^T(\phi_{m_0-2})\geq1-\varepsilon/4\Vert F\Vert_\infty.
\endaligned\end{equation}
By \eqref{eq2.4}, 
there exists a partition $\P_0\in \T^T$ such that for every $\P\geq \P_0$, $I^\P_3(m_0)\leq \varepsilon/3$. 
Let $A_{m}:=\{\gamma\in W^T_x(M):\gamma\subset B_x(m)\}$. Then we have
\begin{equation}\label{eq2.7}\aligned
I^\P_1(m_0)&=\int_{H^T_{x,\P}(M)}(1-\phi_m)F(\gamma)\d\pp^{T,1}_{x,\P}(\gamma)\\
&\leq \Vert F\Vert_\infty\pp^{T,1}_{x,\P}(A^c_{m_0})=\Vert F\Vert_\infty[1-\pp^{T,1}_{x,\P}(A_{m_0-1})]\\
&\leq \Vert F\Vert_\infty[1-\pp^{T,1}_{x,\P}(\phi_{m_0-2})].
\endaligned\end{equation}
In addition, \eqref{eq2.4} implies that
 \begin{equation}\label{eq2.8}\aligned
\lim_{|\P|\rightarrow0}\pp^{T,1}_{x,\P}(\phi_{m_0-2})&=\pp_x^T(\phi_{m_0-2}).
\endaligned\end{equation}
Combining \eqref{eq2.6},\eqref{eq2.7} and \eqref{eq2.8}, there exists a partition $\P_0\subset\P_1\in \T^T$ such that for every $\P\geq \P_1$
\begin{equation}\label{eq2.9}\aligned
I^\P_1(m_0)\leq \varepsilon/3.
\endaligned\end{equation}
Thus, by \eqref{eq2.5}, for the above $\varepsilon$ and or every $\P\geq \P_1$ we have
$$\left|\int_{H^T_{x,\P}(M)}F(\gamma)\d \pp^{T,1}_{x,\P}(\gamma)-\int_{W^T_x(M)}F(\gamma)\d \pp_x^T(\gamma)\right|\leq\varepsilon,$$
which implies \eqref{eq2.13}.
$\hfill\square$

\subsection{Integration by parts formula}
In this subsection, based on Theorem \eqref{T1.1}, we will provide the formula of integration by parts formula for the general non-compact Riemannian manifold. 

Denote by Cameron-Martin space
$$\mathbb{H}=\left\{h\in C([0,T];\mathbb{R}^d): h_0=0,
\|h\|^2_{\mathbb{H}}:=\int_0^T|h_s'|^2\d s<\infty\right\},$$ which is
a separable Hilbert space under $\<h^1,h^2\>_\H:=\int_0^T \langle h'^1_s, h'^2_s\rangle\d s,\
h^1,h^2\in \mathbb H.$ 
For any $h\in \H$, let $\hat{h}$ be the solution of the integrable equation 
\begin{equation}\label{eq2.12}
\hat h_t=h_t +{1\over 2}\Ric_{U_t}\int_0^t h_s\,ds.
\end{equation}
And let $X^{\hat h}(\gamma)\in T_\gamma H^T_x(M)$ be given by
$$X_s^{\hat h}(\gamma):=U_s(\gamma){\hat h}_s(\gamma)\quad \textrm{ for all } s\in [0,T].$$
The vector-valued operator is $\Ric_{u_t}:
\mathbb{R}^d\to\mathbb{R}^d$ is defined by
$$\<(\Ric_{u_t})(a),b\>:=\Ric(u_ta,u_tb),\quad a,b\in\mathbb{R}^d.$$

\beg{thm}\label{T2.2}(Integration by parts formula)Assume that 
\begin{equation}\label{eq2.11}
\mathbb{E}\int_0^T\|\Ric_{u_t}\|^2\d
s<\infty
\end{equation}
where $\|\cdot\|$ is the operator norm from
$\mathbb{R}^d$ to $\mathbb{R}^d$.
 Then for any $h\in \H$, we have
\begin{equation}\label{eq2.12}\aligned
\int_{W^T_x(M)}X^{\hat h}F\d \pp_x^T
=\int_{W^T_x(M)}F\int^T_0
\<h'(t), \d W_t\>\d \pp_x^T
\endaligned\end{equation}
for any $F\in \F C_0^\infty.$
\end{thm}

\beg{proof} We only show that \eqref{eq2.12} holds for all $h\in \H^1$. Otherwise, we may choose a sequence of $h\in \H^1$ such that $h_n\rightarrow h$ in $\Vert\cdot\Vert^{1.2}$. Then the conclusion follows from the dominated convergence theorem. Similar to the argument of the proof of Theorem 1.1. For reader's convenience, we also give a complete proof. 
Given $F\in \F C_0^T$. 
Let $W_x(M_m), H_\P(M_m),\pp^{m,0}_{x,\P}, \pp^{m,1}_{x,\P}$ and $F_m$ be defined as in the proof of Theorem \ref{T1.1}. Thus we have
\begin{equation}\label{eq2.13}\aligned
\int_{W^T_x(M)}X^{\hat h}F_m\d \pp_x^T&=\int_{W_x(M_m)}X^{\hat h}F_m\d \pp^m_x
=\int_{W_x(M_m)}F\int^T_0
\<h'(t), \d W_t\>\d \pp^m_x\\&=\int_{W^T_x(M)}F\int^T_0
\<h'(t), \d W_t\>\d \pp_x^T
\endaligned\end{equation}
where the first equality and the third equality is due to 
$$F_m(\gamma)=0, \gamma \in A^c_m~\text{and}~ \pp_x^T\Big|_{A_m}=\pp_x^{T,m}\Big|_{A_m}.$$ Here $A_m=:\{\gamma\in W^T_x(M); \gamma \subset B_x(m)\}.$ The second equality comes from Theorem 7.16 in \cite{AD}.

Finally, by the dominated convergence theorem, we will obtain \eqref{eq1.12} by letting $R\rightarrow\infty$ in two sides of the equation  \eqref{eq1.13}. \end{proof}


\subsection{An extension to free path spaces}
Let $\mu$ be a probability measure on $M$ and $\pp_\mu$ be the distribution of the Brownian motion starting
from $\mu$ up to time $T$, which is then a probability measure on the free path space
$W^T(M)$.
In fact, we know that
$$\d \pp^T_\mu=\int_M\pp_x^T\d \mu(x),$$
where $\pp_x^T$ is the law of Brownian motion starting at $x$.

Define
\begin{equation*}
\aligned
H^T(M):=\left\{\gamma\in W^T(M): E(\gamma)<\infty\right\}
\endaligned\end{equation*}
and for every $\P\in \T^T$, define
\begin{equation*}
\aligned
H^T_{\P}(M):=\left\{\gamma\in H^T(M)\cap C^2([0,T]/ \P): \nabla \gamma'(t)/ \d t=0, s\notin \P\right\}.
\endaligned\end{equation*}
Setting
\begin{equation*}
\aligned
&\pp^{T,0}_{\mu,\P}:=\int_M\pp^{T,0}_{x,\P}\d \mu(x)\\
&\pp^{T,1}_{\mu,\P}:=\int_M\pp^{T,1}_{x,\P}\d \mu(x).
\endaligned\end{equation*}
Then $\pp^{T,0}_{\mu,\P}$ and $\pp^{T,1}_{\mu,\P}$ are two probability measures on $H^T_{\P}(M)$, so are two probability measures on $W^T(M)$.

\beg{thm}\label{T2.3} Assume that $M$ is complete and stochastic complete. Then 

$(1)$ For each bounded continuous function $F$, we have
\begin{equation}\label{eq2.14}
\aligned
\lim_{|\P|\rightarrow0}\int_{H^T_{\P}(M)}F(\gamma)\d \pp^{T,1}_{\mu,\P}(\gamma)
=\int_{W^T(M)}F(\gamma)\d \pp^T_\mu(\gamma).
\endaligned\end{equation}

$(2)$ 
Suppose that
$$\int_{W^T(M)}\e^{-\frac{1}{6}\int^T_0\Scal(\gamma(t))\d t}\d \pp^T_\mu(\gamma)<\infty.$$
Then for each $F\in \F C_0^\infty$ on $W^T_x(M)$, we have
\begin{equation}\label{eq2.15}
\aligned
\lim_{|\P|\rightarrow0}\int_{H^T_{\P}(M)}F(\gamma)\d \pp^{T,0}_{\mu,\P}(\gamma)
=\int_{W^T(M)}F(\gamma)\e^{-\frac{1}{6}\int^T_0\Scal(\gamma(t))\d t}\d \pp^T_\mu(\gamma),
\endaligned\end{equation}
where $\Scal$ is the scalar curvature of $M$.
\end{thm}

\beg{proof}We only show \eqref{eq2.14}, \eqref{eq2.15} may be derived similarly.
Given a bounded continuous $F$ on $W^T(M)$. By using Theorem \ref{T1.1} and the dominated convergence theorem, we get
\begin{equation}\label{eq2.16}
\aligned
\lim_{|\P|\rightarrow0}\int_{H^T_{\P}(M)}F(\gamma)\d \pp^{T,1}_{\mu,\P}(\gamma)
&=\lim_{|\P|\rightarrow0}\int_M\int_{H^T_{x,\P}(M)}F(\gamma)\pp^{T,1}_{x,\P}(\gamma)\d \mu(x)\\
&=\int_M\lim_{|\P|\rightarrow0}\int_{H^T_{x,\P}(M)}F(\gamma)\pp^{T,1}_{x,\P}(\gamma)\d \mu(x)\\
&=\int_M\int_{W^T_x(M)}F(\gamma)\pp_x^T(\gamma)\d \mu(x)\\
&=\int_{W^T(M)}F(\gamma)\d \pp^T_\mu(\gamma)
\endaligned\end{equation}
\end{proof}

\section{Approximation to Wiener measure for infinite interval}
In this section, we will consider the approximation to Wiener measure for half line. Firstly, we need to introduce some notation.
For a fixed $x\in M$, denote by the path space 
$$W^\infty_x(M):=\{\gamma\in C([0,\infty);M):\gamma(0)=x\}.$$
Then $W^\infty_x(M)$ is a Polish (separable metric)
space under the following uniform distance
$$d_\infty(\gamma,\sigma):=
\displaystyle\sup_{t\in [0,\infty)}\Big(\rho(\gamma(t),\sigma(t))\wedge 1\Big),\quad\gamma,\sigma\in W^\infty_x(M).$$

Let us define $\F C_0^\infty$ the class of cylindric smooth functions on the path space $W^\infty_x(M)$ defined as
\begin{equation}\label{eq3.1}\aligned \F C_0^\infty:=\Big\{F(\gamma)&=f(\gamma(t_1),\cdots,\gamma(t_N)):\ N\geq1,~\gamma\in W^\infty_x(M),\\&
~~~~~~~~~~~~~~~~~~~
0<t_1<t_2\cdots<t_N<\infty,~f\in C_0^{Lip}(M^N)\Big\},\endaligned\end{equation}
where $C_0^{Lip}(M^N)$ is the set of all Lipschitz continuous functions with compact supports on $M^N$. Denote by $\pp_x^\infty$ the Wiener measure on $W^\infty_x(M)$, that is,
 \begin{equation}\label{e3.2}\aligned \int_{W^\infty_x(M)}F(\gamma)\d \pp_x^\infty(\gamma)=\int_{M^N}f(y_1,\cdots,y_N)\prod_i^Np_{\Delta_i}(x_{i-1},x_i)\d x_1\cdots\d x_N,\endaligned\end{equation}
for any $F\in\F C_b^\infty$ of the form $F(\gamma)=f(\gamma(t_1),\cdots,\gamma(t_N))$, where $p_t(x,y)$ is the heat kernel on $M$ and $\Delta_i=t_i-t_{i-1}$.

For any $N\in \N$, let $\P_N=\{0,\frac{1}{2^N},\cdots,\frac{N2^N-1}{2^N},N\}$ and 
\begin{equation}\label{eq3.3}\aligned
&\pp^{1}_{x,N}:=\pp^{N,1}_{x,\P_N}:=\frac{1}{Z^1_{\P_N,N}}\e^{-E/2}\Vol_{G^{N,1}_{x,\P_N}}\\
&\pp^{0}_{x,N}:=\pp^{N,0}_{x,\P_N}=\frac{1}{Z^0_{\P_N,N}}\e^{-E/2}\Vol_{G^{N,0}_{x,\P_N}},
\endaligned\end{equation}
where $Z^0_{\P_N,N}$ and $Z^1_{\P_N,N}$ are normalization constants given by
\begin{equation}\label{eq3.4}\aligned
&Z^0_{\P_N,N}:=\prod^N_{i=1}\left(\frac{\sqrt{2\pi}}{2^N}\right)^{nN2^N}\\
&Z^1_{\P_N,N}:=(\sqrt{2\pi})^{nN2^N}.
\endaligned\end{equation}

Next, the map $\phi_T:H^T_{\P}(M)\mapsto H^\infty_{\P}(M)$ is given by
$$\phi_T(\gamma)=\begin{cases}
&\gamma(t), \quad \quad~ \quad \quad \quad \quad \quad \quad \quad \text{if}~t\in[0,T],\\
&\eta_{\gamma(T)}(t),  \quad \quad \quad \quad \quad \quad \quad\quad \text{if}~t\geq T,
\end{cases}$$
where $\eta_{\gamma(T)}$ is the geodesic with starting point $\gamma(T)$ and initial vector $-\dot{\gamma}(T^-)$.

\beg{thm}\label{T3.1} Assume that $M$ is complete and stochastic complete. Then 

$(1)$ For each bounded continuous function $F$, we have
\begin{equation}\label{eq3.5}\aligned
\lim_{N\rightarrow \infty}\int_{H^{N,x}_{\P_N}(M)}F\circ \phi_N(\gamma)\d \pp^{1}_{x,N}(\gamma)
=\int_{W^\infty_x(M)}F(\gamma)\d \pp_x^\infty(\gamma)
\endaligned\end{equation}

$(2)$ 
Suppose that 
$$\int_{W^T_x(M)}\e^{-\frac{1}{6}\int^T_0\Scal(\gamma(t))\d t}\d \pp_x^T(\gamma)<\infty.$$
Then for  each bounded continuous function $F$, we have
\begin{equation*}\aligned
\lim_{N\rightarrow \infty}\int_{H^{N,x}_{\P_N}(M)}F\circ \phi_N(\gamma)\d \pp^{0}_{x,N}(\gamma)
=\int_{W^T_x(M)}F(\gamma)\e^{-\frac{1}{6}\int^T_0\Scal(\gamma(t))\d t}\d \pp_x^T(\gamma).
\endaligned\end{equation*}
In particular, we have
\begin{equation}\label{eq3.6}\aligned
\lim_{N\rightarrow \infty}\int_{H^{N,x}_{\P_N}(M)}F\circ \phi_N(\gamma)\d \pp^{0}_{x,N}(\gamma)
=\begin{cases}
&0, \quad \quad~ \quad \quad \quad \quad \quad \quad \quad \text{if}~\Scal_M>C,\\
&\int_{W^\infty_x(M)}F(\gamma)\d \pp_x^\infty(\gamma),  \quad \text{if}~\Scal_M\equiv0,\\
&\infty, \quad\quad\quad\quad\quad ~~\quad \quad\quad \text{if}~\Scal_M<-C
\end{cases}
\endaligned\end{equation}
for some constant $C>0$,
where we assume that $0\times\infty=\infty, \infty\pm\infty=\infty$.
\end{thm}

\begin{proof}Here we only show that \eqref{eq3.6} holds(\eqref{eq3.5} may be handled similarly). Let $C_b(W^\infty_x(M))$ be the set of all bounded continuous functions. Since $\F C_0^\infty$ is dense in $C_b(W^\infty_x(M))$, it suffices to prove \eqref{eq3.6} holds for each $F\in \F C_0^\infty$.

We only consider the case of $\Scal_M\equiv0$, the other cases may be proved similarly. Assume that $F\in \F C_0^\infty$ with $F(\gamma)=f(\gamma(t_1),\cdots,\gamma(t_m))$. There exists a $N_0\in \N$ such that $N_0\geq t_m$. Thus $F\in \F C_0^{N_0}$. By (2) of Theorem \eqref{T1.1}, for any $T\geq N_0$ we have
$$\lim_{|\P|\rightarrow0}\int_{H^T_{x,\P}(M)}F(\gamma)\circ \phi_N(\gamma)\d \pp^{T,0}_{x,\P}(\gamma)
=\int_{W^T_x(M)}F(\gamma)\circ\phi_N(\gamma)\d \pp_x^T(\gamma).$$
In particular, we take $T=N_0,N_0+1,\cdots$. By using the form of 
$F(\gamma)=f(\gamma(t_1),\cdots,\gamma(t_m))$ and the definitions of Wiener measures $\pp_x^N$ and $\pp_x^{\infty}$, for each $N\geq N_0$ we have
\begin{equation}\label{eq3.7}\aligned\int_{W^\infty_x(M)}F(\gamma)\d \pp_x^{\infty}(\gamma)&=\int_{W^N_x(M)}F(\gamma)\circ\phi_N(\gamma)\d \pp_x^N(\gamma)\\
&=\lim_{|\P|\rightarrow0}\int_{H^N_{x,\P}(M)}F(\gamma)\circ\phi_N(\gamma)\d \pp^{N,0}_{x,\P}(\gamma).\endaligned\end{equation}
In addition, for each $N$ taking $\P_N=\left\{0,\frac{1}{2^N},\cdots,\frac{N2^N-1}{2^N},N\right\}$. Therefore, by using the law of diagonal lines and \eqref{eq3.7}, we complete the proof.
\end{proof}

\begin{Remark}\label{r3.2} From (2) of Theorem \eqref{T1.1}, we note that the limit to the approximation about $G^{T,1}$ depends on $\int^T_0\Scal(\gamma(t))\d t$. Thus, when $T$ tends to $\infty$ this limit can be explosive or equals to $0$, which also be observed by $(2)$ of Theorem \eqref{T3.1}. Even when $M$ is sphere $\ss^n$ or Hyperbolic space $\H^n$, it can not really characterize the essential geometric property. For this case, we guess that this is because we need to take a suitable new metric on piecewise smooth path space $H_\P^\infty(M)$. In further plan we will continue to work on it. 
\end{Remark}

Next, similar to the argument in subsection 2.3, we also obtain the same results on free path space $W^\infty(M)$, where
$$W^\infty(M)=C([0,\infty);M).$$ Let $\mu$ be a probability measure on $M$ and $\pp_\mu$ be the distribution of the Brownian motion starting
from $\mu$, which is then a probability measure on $W^\infty(M)$.
In fact, we know that
$$\d \pp_\mu=\int_M\pp_x^\infty\d \mu(x).$$

For any $N\in \N$, let $\P_N=\{0,\frac{1}{2^N},\cdots,\frac{N2^N-1}{2^N},N\}$ and
\begin{equation*}
\aligned
&\pp^0_N:=\int_M\pp^{0,x}_{\P_N,N}\d \mu(x)\\
&\pp^1_N:=\int_M\pp^{0,x}_{\P_N,N}\d \mu(x).
\endaligned\end{equation*}
Then $\pp^0_N$ and $\pp^1_N$ are two probability measures on $H^N_{\P_N}(M)$, so are two probability measures on $W^N(M)$. Similar to Theorem \ref{T2.3}, it is easy to obtain the following Theorem \ref{T3.3}.

\beg{thm}\label{T3.3} Assume that $M$ is complete and stochastic complete. Then 

$(1)$ For each bounded continuous function $F$, we have
\begin{equation}\label{eq3.8}
\aligned
\lim_{N\rightarrow \infty}\int_{H^N_{\P_N}(M)}F\circ \phi(\gamma)\d \pp^1_N(\gamma)
=\int_{W^\infty(M)}F(\gamma)\d \pp^\infty_\mu(\gamma).
\endaligned\end{equation}

$(2)$ 
Suppose that 
$$\int_{W^T_x(M)}\e^{-\frac{1}{6}\int^T_0\Scal(\gamma(t))\d t}\d \pp_\mu^T(\gamma)<\infty.$$
Then for  each bounded continuous function $F$, we have
\begin{equation*}\aligned
\lim_{N\rightarrow \infty}\int_{H^N_{\P_N}(M)}F\circ \phi(\gamma)\d \pp^0_N(\gamma)
=\int_{W^T_x(M)}F(\gamma)\e^{-\frac{1}{6}\int^T_0\Scal(\gamma(t))\d t}\d \pp_\mu^T(\gamma).
\endaligned\end{equation*}
In particular,
\begin{equation}\label{eq3.9}\aligned
\lim_{N\rightarrow \infty}\int_{H^N_{\P_N}(M)}F\circ \phi(\gamma)\d \pp^0_N(\gamma)
=\begin{cases}
&0, \quad \quad~ \quad \quad \quad \quad \quad \quad \quad \text{if}~\Scal_M>C,\\
&\int_{W^\infty(M)}F(\gamma)\d \pp_\mu^\infty(\gamma),  \quad \text{if}~\Scal_M\equiv0,\\
&\infty, \quad\quad\quad\quad\quad ~~\quad \quad\quad \text{if}~\Scal_M<-C
\end{cases}
\endaligned\end{equation}
for some constant $C>0$,
where we assume that $0\times\infty=\infty, \infty\pm\infty=\infty$.
\end{thm}

\section*{Acknowledgments}

\beg{thebibliography}{99}

\leftskip=-2mm
\parskip=-1mm

\bibitem{AD} L. Andersson, B. Driver,  \emph{Finite Dimensional Approximations to Wiener Measure and Path Integral Formulas on Manifolds}, J. Funct.Anal. 165(1999), 430--498

\bibitem{A} M. F. Atiyah,  \emph{Circular symmetry and stationary-phase approximation,} Ast ?erisque (1985), no. 131, 43?59, Colloquium in honor of Laurent Schwartz, Vol. 1 (Palaiseau, 1983). MR 816738 (87h:58206) 1.1

\bibitem{CNZ2}V. Cachia, H, Neidhardt,  V. A. Zagrebnov.,\emph{ Comments on the Trotter product formula error-bound estimates for nonself-adjoint semi- groups},Integral Equations Operator Theory, 42(4):425--448, 2002.

\bibitem{CNZ1}V. Cachia, H, Neidhardt,  V. A. Zagrebnov.,\emph{ Accretive perturbations and error estimates for the Trotter product formula},Integral Equations Operator Theory, 39(4):396--412, 2001.

\bibitem{CLW1} X. Chen, X. M. Li, B. Wu, \emph{Functional inequality on loop space over a general Riemannian manifold,}  Preprint.

\bibitem{D} B. Driver,  \emph{A Cameron-Martin type quasi-invariant theorem for Brownian motion on a compact Riemannian manifold},  J. Funct.Anal. 110(1992), 272--376.

\bibitem{F} R. P. Feynman, \emph{Space-Time Approach to Non-Relativistic Quantum Mechanics,} Rev. Mod. Phys.  20(1984), 367-- 387.

\bibitem{FH} R. P. Feynman, Albert R. Hibbs, \emph{Quantum Mechanics and Path Integrals,} emended ed., Dover Publications
Inc., Mineola, NY, 2010, emended and with a preface by Daniel F. Styer.

\bibitem{GJ}L. Glimm, A. Jaffe, \emph{Quantum physics},Springer-Verlag, New York, 1981. A functional integral point of view.

\bibitem{Gr} L. Gross, \emph{Lattice gauge theory; heuristics and convergence. In Stochastic processes¡Xmathematics and physics}, (Bielefeld, 1984), volume 1158 of Lecture Notes in Math., pages 130--140. Springer, Berlin, 1986.

\bibitem{I} W. Ichinose, \emph{On the formulation of the Feynman path integral through broken line paths,}
Comm. Math. Phys. 189 (1997), no. 1, 17--33. MR 1478529 1.1

\bibitem{LZH} Z. H. Li, \emph{A finite dimensional approximation to pinned Wiener measure on some symmetric spaces,} {\it arXiv:1702.06747v1}.

\bibitem{AL} Adrian P. C. Lim, \emph{Path integrals on a compact manifold with non-negative curvature,} Rev. Math. Phys. 19 (2007), no. 9, 967?1044. MR 2355569 (2008m:58080) 1.1

\bibitem{NZ} H. Neidhardt, V. A. Zagrebnov., \emph{Trotter-Kato product formula and operator-norm convergence}, Comm. Math. Phys., 205(1)(1999):129--159, .

\bibitem{P} Mark A. Pinsky,  \emph{Isotropic transport process on a Riemannian manifold}, Trans. Amer. Math. Soc. 218 (1976), 353--360. MR 0402957 1.1

\bibitem{TL} Thomas Laetsch, \emph{An approximation to Wiener measure and quantization of the Hamiltonian on manifolds with non-positive sectional curvature,} J. Funct. Anal. 265 (2013), no. 8, 1667--1727. MR 3079232 1.1, A

\end{thebibliography}
\end{document}